# Stochastic dynamics of the triple-well potential systems driven by colored noise

Yanxia Zhang, Yanfei Jin[*]

*Department of Mechanics, Beijing Institute of Technology, Beijing 100081, China*

**Abstract:** A stochastic averaging technique based on energy-dependent frequency is extended to dynamical systems with triple-well potential driven by colored noise. The key procedure is the derivation of energy-dependent frequency according to the four different motion patterns in triple-well potential. Combined with the stochastic averaging of energy envelope, the analytical stationary probability density (SPD) of tri-stable systems can be obtained. Two cases of strongly nonlinear triple-well potential systems are presented to explore the effects of colored noise and validate the effectiveness of the proposed method. Results show that the proposed method is well verified by numerical simulations, and has significant advantages, such as high accuracy, small limitation and easy application in multi-stable systems, compared with the traditional stochastic averaging method. Colored noise plays a significant constructive role in modulating transition strength, stochastic fluctuation range and symmetry of triple-well potential. While, the additive and multiplicative colored noises display quite different effects on the features of coherence resonance (CR). Choosing a moderate additive noise intensity can induce CR, but the multiplicative colored noise cannot.

**I. INTRODUCTION**

Numerous physical systems are associated with the random fluctuating environment or noise, which can induce many complicated nonlinear phenomena, such as stochastic transition response and resonance phenomenon, even in simple systems. Specially, coherence resonance (CR), called as stochastic resonance (SR) without external periodic excitation initially, describes the phenomenon that the output response in an excitable system can become quite regular with the addition of a moderate amount of noise. Thus, the study of dynamical systems under random excitations has attracted accelerating attention because of its wide applications in many fields [1-5] and it is necessary to consider the effects of random noise excitations on the stochastic dynamical behaviors of the systems. Recently, the stochastic response and CR of the dynamical systems with random excitations have been studied in detail by many researchers [4-11]. For example, Xu et al. [4] found that the additive noise can enlarge the transition response in the studying of a bi-stable system under white noise. Xiao et al. [5] studied the effect of the correlated white noise on the response of the mono-stable energy harvesting system, and found that correlated noise can improve the energy harvesting capacity. Bogatenko et al. [9] discussed the CR in an excitable potential well by the addition of broadband noise. Jin and Xu [10,11] studied the noise-induced response and CR in a one-dimensional dynamical system driven by white noise, and found that CR occurs by choosing an optimal noise intensity. However, most of the previous studies focus on mono- and bi-stable systems under white noise, less attention is paid to the stochastic response and CR phenomenon induced by colored noise in multi-stable dynamical systems. To explore the stochastic dynamics of the multi-stable dynamical systems driven by colored noise, it is important to develop analytical approaches for solving the stationary probability density.

Currently, an effective and powerful method in studying stochastic dynamical systems is the stochastic averaging method, which is proposed initially by Stratonovich [12] and is validated rigorous by Khasminskii [13] and Papanicolaou [14]. Two approximation procedures are carried out in the method. One is the approximation of the excitation process by Gaussian white noise to make the system response described by a Markov process, whose probability density is governed by Fokker-Planck-Kolmogorov (FPK) equation. Another is the time averaging to eliminate the fast varying processes and reduce the system dimension. The stochastic averaging method has been used and developed by many researchers [15-18] for the single-degree-of-freedom and the multi-degree-of-freedom systems. For example, Spanos et al. [16,17] not only developed and applied the stochastic averaging method in various nonlinear random vibrations but also demonstrated the merit of the approximate solution by deriving the probability density function of the system response. Zhu and Lin [18] proposed a stochastic averaging called the stochastic averaging of energy envelope by considering the conservative part of the Wong-Zakai correction terms of system. Zhu et al. [19] extended the stochastic averaging of energy envelope method and applied into Hamilton systems. Furthermore, Zhu and Huang [6] proposed an averaging method to deal with strongly nonlinear stochastic systems subjected to random noise excitations. Using the generalized harmonic functions, they applied the proposed stochastic averaging method to analyze the stationary response of a typical Duffing-van der Pol oscillator under wide-band random excitation, and verified the accuracy of the new method, which depends on the bandwidth of the excitation, by numerical simulation. Whereas, the method is suitable for broadband excitation but not for narrow-band excitation, such as colored noise.

---

[*] Corresponding author.
 E-mail address: 3120170019@bit.edu.cn (Y. Zhang), jinyf@bit.edu.cn (Y. Jin).

Colored noise with non-zero correlation time can better describe the actual environmental random excitations especially when the correlation time of random fluctuation is large. Recently, theoretical studies for stochastic dynamical systems with colored noise are mainly focused on the mono- and bi-stable systems [20-24]. Xu et al. [20] investigated the dynamical behavior of a bi-stable Duffing–van der Pol oscillator with colored noise by using the stochastic averaging method and revealed some interesting nonlinear phenomena through the SPD of amplitude obtained by FPK equation. Liu et al. [23] studied the mono-stable response of a nonlinear vibration energy harvester under colored noise by using the extended stochastic averaging of energy envelope and found that colored noise has a significant effect on the performance of energy harvesting. In addition, Zhu et al. [24] proposed a procedure to apply the stochastic averaging method to systems with double-well potential by using the correlation functions of the excitation processes instead of the power spectral densities and adopting the variable natural period and frequency of system corresponding to different energy levels. Xu et al. [4] extended the stochastic averaging to the bi-stable nonlinear potential energy with the energy-dependent frequency and validated the precision of the method by Monte Carlo simulations (MCS). However, compared with mono- or bi-stable system, the study of dynamical behaviors in tri-stable systems is more important because such systems have a wide range of applications in various fields, such as engineering [25,26], optomechanical system [27] and condensed matter physics [28]. Yet, since the strongly nonlinear potential function considered in a tri-stable system has the triple-well shape, the system motion is more complicated. The total energy, the initial condition and the shape of potential function determine the forms of motion, i.e., the particle moves in one of three potential wells or jumps among three wells. Moreover, the consideration of colored noise in a strongly nonlinear tri-stable system brings the complexity in the mathematical calculation and makes the traditional stochastic averaging methods invalid.

Inspired by this issue, this paper is devoted to developing an improved stochastic averaging method to study the effects of colored noise on the stochastic dynamics of the tri-stable systems. The rest of paper is organized as follows. Sec. II derives the energy-dependent period and the energy-dependent frequency according to the four different motion patterns in the triple-well potential, and calculates the drift and diffusion coefficients of the energy process. Then, the SPD can be obtained to analyze the dynamics of a tri-stable system under colored noise. In Sec. III, two illustrative cases are given to explore the effects of colored noise and validate the precision of the proposed method. Sec. IV explores the effects of colored noise on the CR by the power spectrum density (PSD) and the quality factor. Finally, Sec. V is devoted to some specific conclusions.

**II. THE STOCHASTIC AVERAGING FOR THE TRIPLE-WELL POTENTIAL WITH COLORED NOISE**

The stochastic system of concern is a strongly nonlinear triple-well potential driven by colored noise excitations. The system is governed by

$$\ddot{X} + \frac{dU(X)}{dX} = h(X, \dot{X}) + g_1(X)\xi_1(t) + g_2(X)\xi_2(t),$$
$$U(X) = \frac{1}{2}k_1 X^2 - \frac{1}{4}k_2 X^4 + \frac{1}{6}k_3 X^6,$$
(1)

where $U(X)$ is the triple-well potential function; the positive parameters $k_1$, $k_2$ and $k_3$ are the linear, cubic and quintic stiffness coefficients, respectively. $h(X, \dot{X})$ denotes the damping force; $g_i(X)$ denotes deterministic function that characterizes the state-dependent action of noise terms $\xi_i(t)$ $(i=1,2)$.

The noise terms $\xi_i(t)$ $(i=1,2)$ considered in this work are focused on colored noise, and their statistical properties are characterized as follows

$$\langle \xi_i(t) \rangle = 0, \quad \langle \xi_i(t)\xi_i(s) \rangle = \frac{D_i}{\tau_i} \exp\left[-\frac{|t-s|}{\tau_i}\right],$$
$$\langle \xi_1(t)\xi_2(s) \rangle = \langle \xi_2(t)\xi_1(s) \rangle = \frac{\lambda\sqrt{D_1 D_2}}{\tau_1 \tau_2} \exp\left[-\frac{|t-s|}{\tau_1 \tau_2}\right],$$
$$S_i(\omega) = \int_{-\infty}^{+\infty} \frac{D_i}{\tau_i} e^{-|\tau'|/\tau_i} e^{-i\omega\tau'} d\tau' = \frac{D_i}{1+\tau_i^2 \omega^2},$$
(2)

where $D_i$ and $\tau_i$ denote the noise intensity and the correlation time of the colored noise, respectively; $\lambda$ is noise cross-correlation; $S_i(\omega)$ is the spectral density of colored noise.

**A. Energy-dependent frequency**

The deterministic conservative system of Eq. (1) is written as

$$\ddot{X} + \frac{dU(X)}{dX} = 0.$$
(3)

The total energy of the system is given by

$$H(X, \dot{X}) = \frac{1}{2}\dot{X}^2 + \frac{1}{2}k_1 X^2 - \frac{1}{4}k_2 X^4 + \frac{1}{6}k_3 X^6,$$
(4)

where the positive parameters $k_1$, $k_2$ and $k_3$ are the linear, cubic and quintic stiffness coefficients, respectively. When $k_2^2 - 4k_1 k_3 < 0$, only one stable equilibrium at origin exists and the potential function is mono-stable. When $k_2^2 - 4k_1 k_3 > 0$ (the case we focused on), as shown in Figs. 1(a) and 1(b), the potential function is tri-stable with two symmetric stable equilibria, one stable equilibrium at origin and two unstable saddle points which are given by

$$X_{s1} = \sqrt{\frac{k_2 + \sqrt{k_2^2 - 4k_1 k_3}}{2k_3}}, X_{s2} = -\sqrt{\frac{k_2 + \sqrt{k_2^2 - 4k_1 k_3}}{2k_3}}, X_{s3} = 0, X_u = \pm\sqrt{\frac{k_2 - \sqrt{k_2^2 - 4k_1 k_3}}{2k_3}}.$$
(5)

Figure 1(c) shows the periodic motions with different energy levels. If the total energy is more than the critical potential energy $U_1$, the periodic motion jumps among the three wells. Otherwise, the periodic motion is restricted in one of three potential wells - the right one, the left one or the middle one depending on the initial condition. It should be noted here that the periodic motion will be only restricted in the deepest potential well if total energy is less than $U_2$. Thus, among these four kinds of motions, the periodic trajectory in which one depends on the total energy, the initial condition and the shape of potential function.

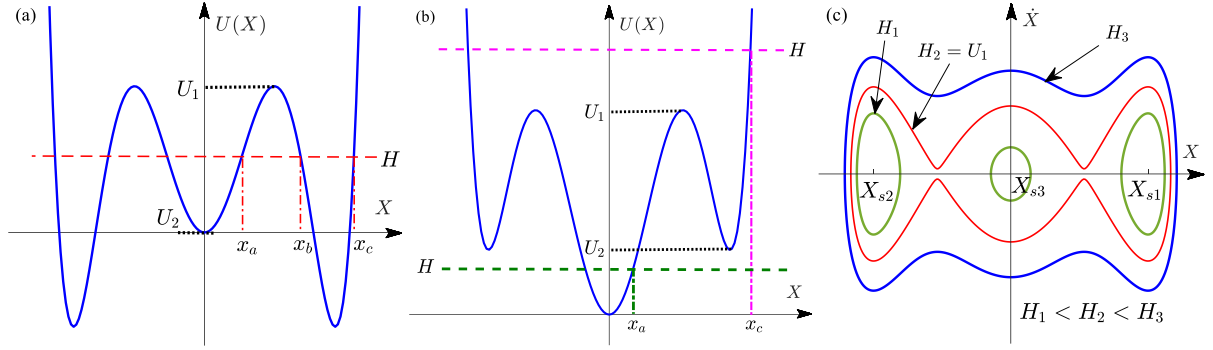

**FIG. 1.** (a) The potential function for fixed $k_1 = 1$, $k_2 = 4.5$ and $k_3 = 3.5$; (b) The potential function for fixed $k_1 = 1$, $k_2 = 4.5$ and $k_3 = 4$; (c) Periodic motions with different energy levels.

For a given energy level $H > U_1$, the energy-dependent period of the motion can be expressed as

$$T(H) = \oint dt = \oint \frac{dX}{\dot{X}} = 4 \int_0^{x_c} \frac{dX}{\sqrt{2H - 2U(X)}}, \quad (6)$$

where $x_c$, along with $x_a$ and $x_b$ as shown in Fig. 1, denotes the min or the max displacement determined by $H = U(A)$.

For a given energy level $U_2 < H < U_1$, the initial condition determines whether the motion is in middle or bilateral potential well. Thus, the energy-dependent period of the motion can be expressed respectively as

$$T(H) = \oint dt = \oint \frac{dX}{\dot{X}} = 4 \int_0^{x_a} \frac{dX}{\sqrt{2H - 2U(X)}}, \quad (7)$$

$$T(H) = \oint dt = \oint \frac{dX}{\dot{X}} = 2 \int_{x_b}^{x_c} \frac{dX}{\sqrt{2H - 2U(X)}}. \quad (8)$$

For a given energy level $H < U_2$, the shape of potential function, whose well is the deepest, determines whether the motion is in middle or bilateral potential well. For the shape in Fig. 1(a), the period expression becomes Eq. (8). Whereas, for the shape in Fig. 1(b), the period expression is given as Eq. (7).

Then, the energy-dependent frequency of the motion is given as

$$\omega(H) = \frac{2\pi}{T(H)}. \quad (9)$$

Figure 2 shows the behavior variation of the total energy $H$ and the energy-dependent frequency $\omega(H)$ for the state variables of tri-stable system. It can be seen that $\omega(H)$ has two symmetric peaks at bilateral stable points and a peak at origin point. That is, in the attraction domain of the stable point, as $H$ increases, $\omega(H)$ decreases. Whereas, out of the attraction domain, $\omega(H)$ turns to an increasing trend suddenly with the increasing energy, which is the reason that the motion jumps from a small trajectory restricted in one well to a large trajectory moving among three wells. Evidently, it can be seen in Fig. 2(b), the energy-dependent frequency of the tri-stable system, proposed in this paper, is very different from those of mono-stable and bi-stable systems [4,23,24].

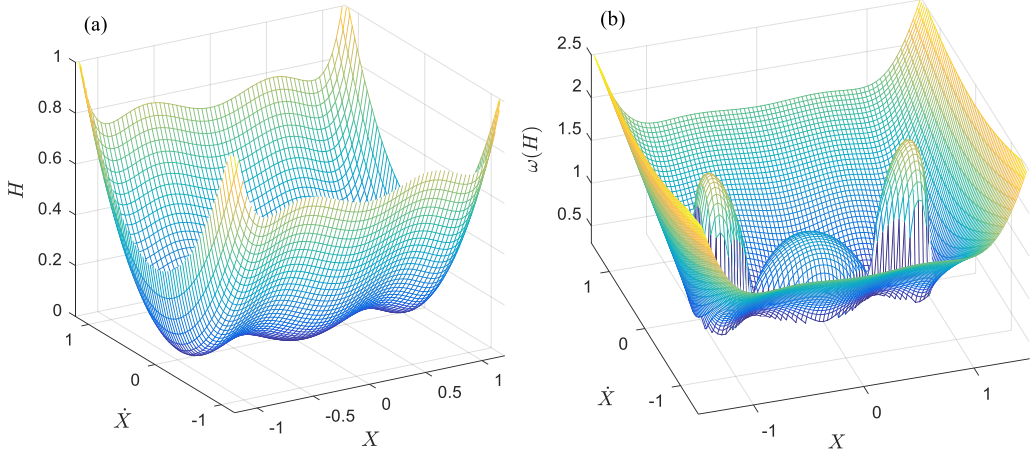

**FIG. 2.** (a) The total energy and (b) the energy-dependent frequency for the state variables of unperturbed system for fixed $k_1 = 1$, $k_2 = 4.5$ and $k_3 = 4$.

### B. Stochastic averaging procedure

Based on the above four different motion patterns, a stochastic averaging technique depends on energy-dependent frequency is proposed to solve stochastic tri-stable dynamical systems. The system (1) can be rewritten as following first-order equations associated with system velocity and total energy [29,30]

$$\dot{X} = \pm\sqrt{2H - 2U(X)}, \tag{10}$$

$$\dot{H} = -\dot{X}h(X,\dot{X}) + \dot{X}\sum_{i=1}^{n} g_i(X,\dot{X})\xi_i(t). \tag{11}$$

Since the energy process $H(t)$ is approximate Markovian, the Itô equation for the limiting diffusion process $H(t)$ can be expressed as the following form

$$dH = m(H)dt + \sigma(H)dB(t), \tag{12}$$

where $B(t)$ is a standard Wiener process, $m(H)$ and $\sigma(H)$ are the drift and diffusion coefficients of the energy process, respectively, which can be calculated by

$$m(H) = -\langle \dot{X} \cdot h(X,\dot{X})\rangle_t + \int_{-\infty}^{0}\sum_{i,j=1}^{2}\left\langle \left[\dot{X}g_j(X)\right]_{t+\tau}\frac{\partial}{\partial H}\left[\dot{X}g_i(X)\right]_t\right\rangle R_{ij}(\tau)d\tau, \tag{13}$$

$$\sigma^2(H) = \int_{-\infty}^{+\infty}\sum_{i,j=1}^{2}\left\langle \left[\dot{X}g_j(X)\right]_{t+\tau}\left[\dot{X}g_i(X)\right]_t\right\rangle R_{ij}(\tau)d\tau, \tag{14}$$

where $R_{ij}(\tau)$ denotes the correlation function of noise excitation; $\langle\cdot\rangle_t = \frac{1}{T}\int_0^T (\cdot)dt$ denotes the time averaging. According to the different total energy levels, the closed-loop integration needs to be carried out corresponding to the different periods of motion developed in Eqs. (6)-(8).

Then, the stationary probability density (SPD) of the total energy can be obtained from Eq. (12) as

$$p(H) = \frac{C_0}{\sigma^2(H)}\exp\left[\int \frac{2m(H)}{\sigma^2(H)}dH\right], \tag{15}$$

where $C_0$ is a normalization constant.

It can be proved that

$$\frac{d}{dH}\ln\left[T(H)\langle \dot{X}^2\rangle_t\right] = \frac{1}{\langle \dot{X}^2\rangle_t}. \tag{16}$$

Thus, the joint SPD of the tri-stable system (1) can be derived as

$$p(X,\dot{X}) = \frac{p(H)}{T(H)} = \frac{C_0 \langle \dot{X}^2 \rangle_t}{\sigma^2(\mu)} \exp\left[\int_0^H \left(\frac{2m(\mu)}{\sigma^2(\mu)} - \frac{1}{\langle \dot{X}^2 \rangle_t}\right) d\mu\right]\Bigg|_{H=\frac{1}{2}\dot{X}^2+U(X)}. \tag{17}$$

From Eq. (17), the marginal SPD are given as

$$p(X) = \int_{-\infty}^{\infty} p(X,\dot{X}) \, d\dot{X},$$
$$p(\dot{X}) = \int_{-\infty}^{\infty} p(X,\dot{X}) \, dX. \tag{18}$$

## III. STOCHASTIC RESPONSE ANALYSIS

To explore the effects of colored noise on the stochastic response of the triple-well potential systems and to validate the effectiveness of the above proposed method, two cases of nonlinear tri-stable potential systems driven by colored noise are presented in this section. The analytical expressions of SPD are derived from Eqs. (17)-(18) and the numerical simulations are obtained through MCS from the original systems. Furthermore, comparisons of the proposed method and the traditional stochastic averaging methods are also discussed in the case I for illustrating its efficiency.

### A. Case I

A class of strongly nonlinear tri-stable potential models have been a hot topic in recent researches and have a wide range of applications in many practical systems, such as energy harvester systems [25,31,32] used for harvesting electric power from ambient environment. A general tri-stable system driven by colored noise is mainly considered firstly, i.e., $g_1(X)=1$, $g_2(X)=0$ and $h(X,\dot{X})=\beta\dot{X}$. Then, system (1) can be rewritten as

$$\ddot{X} + \beta\dot{X} + k_1 X - k_2 X^3 + k_3 X^5 = \xi_1(t). \tag{19}$$

Here, Eq. (19) can be treated as a non-dimensionalized equation of a general inductive-type energy generator subjected to colored noise excitation [31].

The drift and diffusion coefficients of the energy process obtained from Eqs. (13)-(14) are given by

$$m(H) = -\beta \langle \dot{X}^2 \rangle_t + \frac{D_1}{1+\tau_1^2 \omega^2(H)}, \tag{20}$$

$$\sigma^2(H) = \frac{2D_1}{1+\tau_1^2 \omega^2(H)} \langle \dot{X}^2 \rangle_t. \tag{21}$$

Substituting Eqs. (20)-(21) into Eq. (15) and Eq. (17), the SPD of total energy and the joint SPD of displacement and velocity can be obtained as

$$p(H) = \frac{C_0(1+\tau_1^2 \omega^2(H))}{D_1 \omega(H)} \exp\left[-\frac{\beta(1+\tau_1^2 \omega^2(H))}{D_1} H\right], \tag{22}$$

$$p(X,\dot{X}) = \frac{C_0(1+\tau_1^2 \omega^2(H))}{D_1} \exp\left[-\frac{\beta(1+\tau_1^2 \omega^2(H))}{D_1}\left(\frac{1}{2}\dot{X}^2 + \frac{1}{2}k_1 X^2 - \frac{1}{4}k_2 X^4 + \frac{1}{6}k_3 X^6\right)\right]. \tag{23}$$

Then, the marginal SPD of system displacement and velocity can be further obtained by substituting Eq. (23) into Eq. (18).

In addition, owing to $p(a) = p(H)\left|\frac{dH}{da}\right|_{H=U(a)}$, the SPD of amplitude can be derived as follows

$$p(a) = \frac{C_0(1+\tau_1^2\omega^2(H))}{D_1\omega(H)}|k_1 a - k_2 a^3 + k_3 a^5|\exp\left[-\frac{\beta(1+\tau_1^2\omega^2(H))}{D_1}(\frac{1}{2}k_1 a^2 - \frac{1}{4}k_2 a^4 + \frac{1}{6}k_3 a^6)\right]. \tag{24}$$

The parameters in this case are chosen as $k_1=1$, $k_2=4.5$, $k_3=4$, $\beta=0.1$, $D_1=0.01$ and $\tau_1=0.5$, unless otherwise mentioned. In Figs. 3 and 4, solid lines denote the analytical results obtained from Eq. (23) by the proposed stochastic averaging, while circle symbol represents the numerical results obtained from Eq. (19) through MCS. It is obvious that the theoretical results are fairly in good agreement with those from numerical simulations. The effects of nonlinear stiffness coefficients on SPD are analyzed in Fig. 3(c) and 3(d). It is seen that the stochastic transition is constantly moving from bilateral stable points toward the origin stable point as $k_3$ increases, and the topological structure of SPD changes from bimodal to tri-modal distribution then end with mono-modal distribution. Whereas, the change is completely opposite as $k_2$ increases. It is the reason that the depth of the bilateral potential wells decreases with increasing $k_3$ but increases with increasing $k_2$. That is, the nonlinear stiffness coefficients can change the topological structure of SPD and induce stochastic P-bifurcation.

Figure 4 presents the effects of colored noise on SPD. As shown in Figs. 4(a) and 4(b), $p(x)$ at two unstable saddle points $x_u \approx \pm 0.53$ increases and the distribution of $p(\dot{x})$ is broadened with increasing $D_1$, which indicates that increasing noise intensity can strengthen the transition among the three potential wells and enlarge the fluctuation range of stochastic response. Whereas, the increase of correlation time $\tau$ weakens both the transition among the three potential wells and the fluctuation range of stochastic response, as shown in Figs. 4(c) and 4(d). In other words, the noise intensity and time correlation play an opposite role on the SPD.

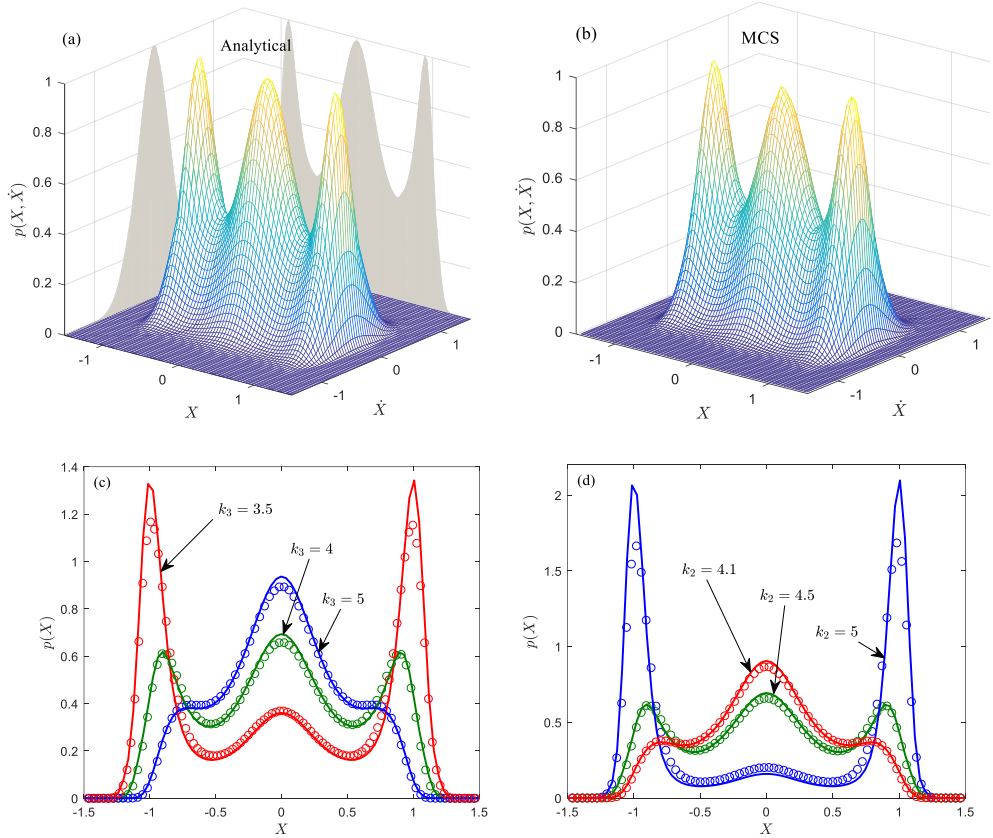

**FIG. 3.** (a) The joint SPD by Eq. (23) from the proposed stochastic averaging; (b) the joint SPD by Eq. (19) from MCS; (c) the marginal SPD for different $k_3$ values; (d) the marginal SPD for different $k_2$ values. Solid line: the proposed stochastic averaging; Circle symbol: MCS.

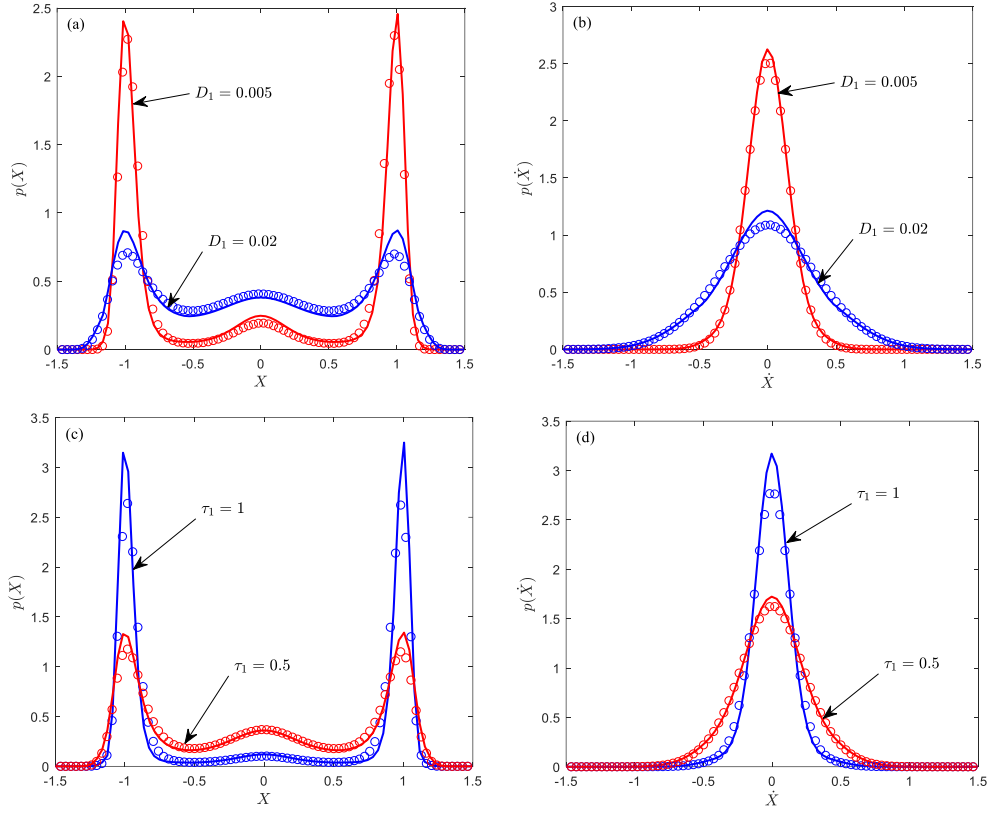

**FIG. 4.** The marginal SPD of (a) system displacement and (b) system velocity for different noise intensities for fixed $\tau = 0.5$ and $k_3 = 3.5$; The marginal SPD of (c) system displacement and (d) system velocity for different correlation times for fixed $D_1 = 0.01$ and $k_3 = 3.5$. Solid line: the proposed stochastic averaging; Circle symbol: MCS.

Comparisons of the proposed method and the traditional stochastic averaging method are presented in Fig. 5. In Fig. 5, the solid line and the dashed line are the analytical results obtained by our proposed method and the stochastic averaging method of amplitude in Ref. [6], respectively. And the circle symbol denotes the numerical results of Eq. (19). Fig. 5(a) shows the SPD of amplitude with parameters $k_2 = 4$, $D_1 = 0.001$ and $\tau_1 = 0.1$. Apparently, the improved method proposed by us shows better accuracy with MCS compared with the method presented in Ref. [6]. Moreover, our proposed method has the advantage of a wider applicability for strongly nonlinear tri-stable systems. For example, we employed the method proposed in Ref. [6] to calculate the averaged FPK equation and obtain the SPD of amplitude of Eq. (19). However, as pointed out in Ref. [31], the valid range of theoretical result obtained by the method in Ref. [6] is very narrow. That is, the condition $k_2^2 - 40k_1k_3/9 < 0$ should be satisfied which greatly reduces the stability region of tri-stable case for system (19). While the theoretical result (24) does not need to meet the above restriction. On the other hand, by using the method proposed in Ref. [6], the SPD of total energy $p(H)$ of system (19) cannot be obtained because of the inexistence of the inverse function $a = U^{-1}(H)$. To avoid this tough problem, the SPD of amplitude (24) can be obtained from the SPD of total energy directly. Fig. 5(b) shows the SPD of total energy whose theoretical and numerical results are in good agreement. Therefore, compared with the traditional method in Ref. [6], our proposed method has the significant advantages, such as high accuracy, small limitations and easy application in tri-stable system.

Furthermore, some recent works, e.g. Refs. [21, 26], assumed that the system oscillates at one frequency and set the frequency as 1 during the procedure of stochastic averaging method. As seen in Fig. 2(b), the frequency is not a fixed value but varies with the state variables of the system. Hence, it is necessary to make a comparison between the proposed method in this paper and the method presented in Ref. [26] where the frequency is set as 1. Fig. 6 shows the marginal SPD of system displacement with three kinds of results, i.e., the numerical results (circle symbol), the theoretical results obtained by the proposed method (solid lines) and by the method in Ref. [26] (dashed lines). One can find in Fig. 6(a) that the solid and dashed lines almost overlap and these two kinds of theoretical results are consistent with numerical results very well for $k_3 = 5$, where the corresponding potential function has the deepest middle well. Whereas, if the bilateral well of potential function is deeper than the middle one (i.e. $k_3 = 3.5$), the particle is mainly confined in the movement between bilateral potential wells, as shown in Fig. 6(b). For this case, the proposed method offered better accuracy than the method in Ref. [26], especially at the position of $x = 0$. Thus, the comparison indicates that the adopting of the energy-dependent frequency makes the theoretical results more accurate and effective, which is independent of the shape of potential function.

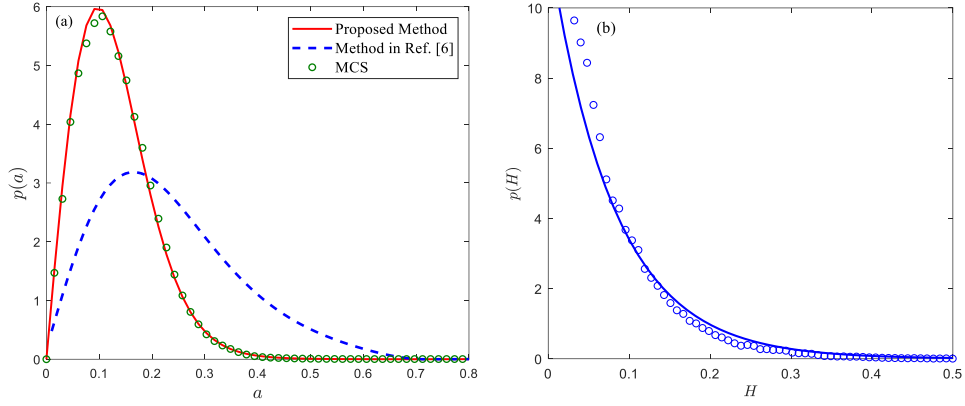

**FIG. 5.** The SPD of (a) amplitude and (b) total energy. Solid line: the analytical solution Eqs. (22) and (24) by using our proposed method; Dashed line: the analytical solution obtained by using the method proposed in Ref. [6]; Circle symbol: MCS.

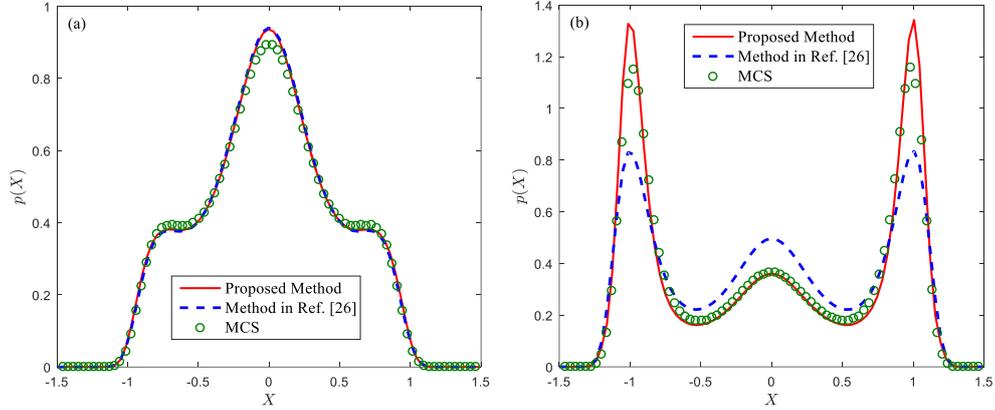

**FIG. 6.** The marginal SPD of system displacement for fixed (a) $k_3 = 5$ and (b) $k_3 = 3.5$.

## B. Case II

Many previous studies assumed that the additive and multiplicative noises have different noise sources, and they are independent of each other. Whereas, in some situations, additive and multiplicative noises may come from the same noise source, which makes them have some kind of noise cross-correlation. Moreover, recent studies [5,33] found that the cross-correlation between the white noises has a great influence on the dynamic behavior of systems. Thus, the effects of the cross-correlation between the additive and multiplicative colored noises on the strongly nonlinear tri-stable system are mainly studied in this case, i.e., $g_1(X) = 1$, $g_2(X) = X$ and

$h(X,\dot{X}) = (\beta + \beta_1 X^2)\dot{X}$. Then, system (1) can be rewritten as

$$\ddot{X} + (\beta + \beta_1 X^2)\dot{X} + k_1 X - k_2 X^3 + k_3 X^5 = \xi_1(t) + X\xi_2(t), \quad (25)$$

where $\beta$ and $\beta_1$ represent the linear and nonlinear damping coefficient, respectively. It should be noted that the colored noise can degenerate into white noise when correlation time tends to zero. Hence, Eq. (25) can be also regarded as a nonlinear system subjected to white noises or combination of white noise and colored noise.

The drift and diffusion coefficients of the energy process for system (25), which can be obtained from Eqs. (13) and (14), are presented as follows

$$m(H) = -\left\langle (\beta + \beta_1 X^2)\dot{X}^2 \right\rangle_t + \frac{D_1}{1+\tau_1^2\omega^2(H)} + \frac{D_2}{1+\tau_2^2\omega^2(H)}\left\langle X^2 \right\rangle_t + \frac{2\lambda\sqrt{D_1 D_2}}{1+\tau_1^2\tau_2^2\omega^2(H)}\left\langle X \right\rangle_t, \quad (26)$$

$$\sigma^2(H) = \frac{2D_1}{1+\tau_1^2\omega^2(H)}\left\langle \dot{X}^2 \right\rangle_t + \frac{2D_2}{1+\tau_2^2\omega^2(H)}\left\langle \dot{X}^2 X^2 \right\rangle_t + \frac{4\lambda\sqrt{D_1 D_2}}{1+\tau_1^2\tau_2^2\omega^2(H)}\left\langle \dot{X}^2 X \right\rangle_t. \quad (27)$$

Then, the joint SPD and marginal SPD of system (25) is derived by substituting Eqs. (26)-(27) into Eqs. (17)-(18).

The effects of cross-correlation between colored noises on the SPD for the tri-stable system are analyzed in this subsection by choosing parameters as $k_1 = 1$, $k_2 = 4.5$, $k_3 = 4$, $\beta = 0.1$, $\beta_1 = 0.05$, $D_1 = 0.005$, $D_2 = 0.005$, $\tau_1 = 0.5$, $\tau_2 = 0.5$. Figure 7 shows the analytical and MCS results of SPD with different noise cross-correlation $\lambda$. In Figure 7(d), the analytical results obtained from the proposed stochastic averaging are marked by solid lines while circle symbol represents the numerical results through MCS. Obviously, the analytical and the numerical results are in good agreement, which well validates the precision of the proposed method. As shown in Figures. 7(c) and 7(d), the symmetry of SPD is broken by $\lambda$, and the curve tilts more severely with increasing $|\lambda|$. Besides, the tilted curves are anti-symmetric for the positive and negative value of $\lambda$. It indicates that noise cross-correlation can break the symmetry of potential wells and induce phase transition. From Figure 7, it shows that the proposed method can be successfully applied to analyze the tri-stable system with various noises, including uncorrelated and correlated additive and multiplicative colored noises and white noises.

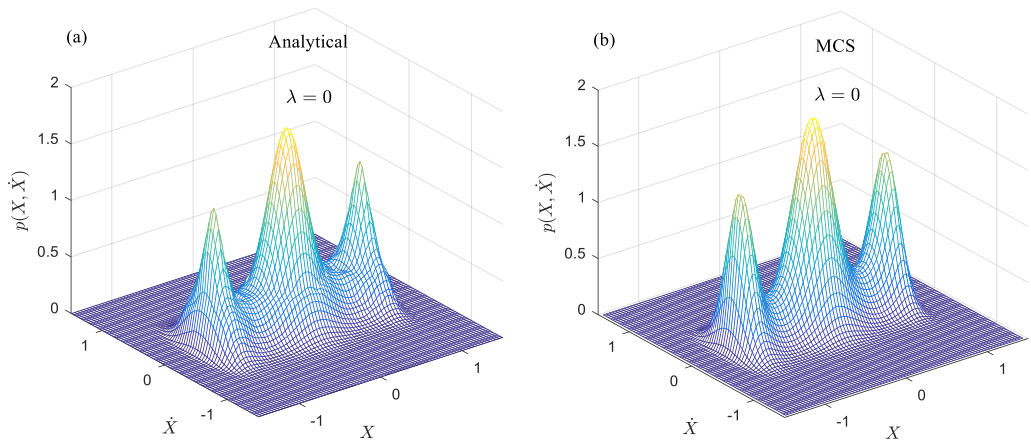

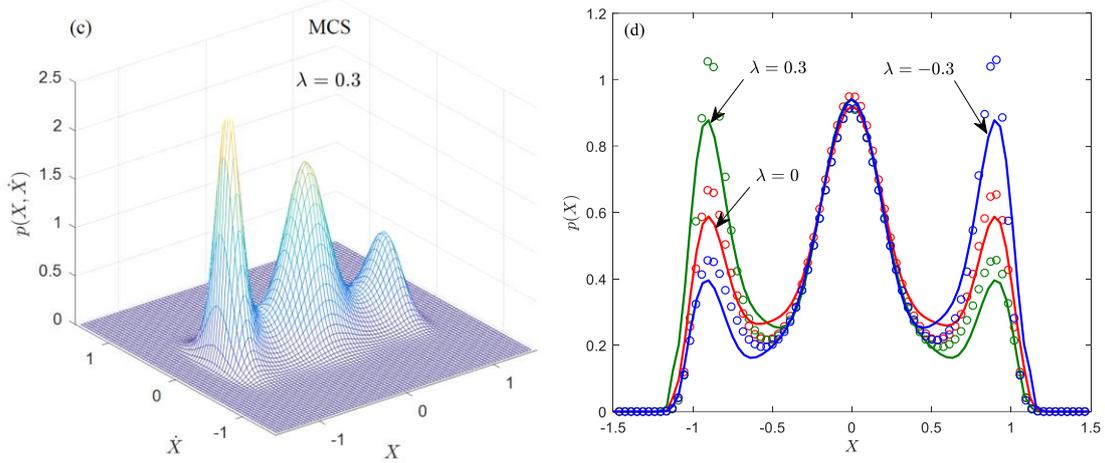

**FIG. 7.** (a)The joint SPD of system (25) obtained by the proposed stochastic averaging; (b) and (c) the joint SPD of system (25) obtained through MCS; (d) The marginal SPD of system (25) for different noise cross-correlations. Solid line: the proposed stochastic averaging; Circle symbol: MCS.

## IV. COHERENCE RESONANCE

To explore the effect of colored noise on the coherence resonance (CR), the CR in the system (1) can be characterized by the power spectrum density (PSD) and the quality factor. The quality factor is defined as $\eta = h\omega_m / \Delta\omega_f$, where $h$ is the maximal peak height of the power spectrum density, $\omega_m$ is the peak frequency and $\Delta\omega_f$ is the width of the power spectrum measured at the height of $h/\sqrt{e}$.[7]

For purpose of exploring the effect of colored noise on the CR, the system (19) in case I is considered here to reveal the distinctive characteristic of CR, as shown in Fig. 8, and the parameters take the same value as Case I. In Fig. 8(a), the PSD as a function of frequency $\omega_f$ is presented with different noise intensity $D_1$. It can be seen that the maximal peak value of PSD increases at first and then decreases with the increase of $D_1$, which indicates that there is an optimal noise intensity to maximize PSD. Obviously, it is a characteristic of CR. For further understanding the variation tendency of the peaks in PSD in Fig. 8(a), the time series of system displacement with the corresponding different $D_1$ are present in Fig. 8(b). It is clear that for a small noise intensity $D_1 = 0.001$, there is only intrawell motion, which explains why there is only one peak for the blue curve in Fig. 8(a). Meanwhile, it indicates that the resonance peak with the blue curve denotes the intrawell resonance. For $D_1 = 0.01$ (see Fig. 8(b)), the particle moves among the three potential wells and it mainly performs interwell resonance (see red curve in Fig. 8(a)). Then, the interwell CR phenomenon is weaken and the intrawell resonance become strong with increasing $D_1$ to 0.02. Figure 8(c) displays the non-monotonic variation of the quality factor $\eta$ with $D_1$, which indicates that CR occurs and the system performs the best coordinated oscillatory character at a moderate additive noise intensity $D_1 \approx 0.01$. Figure 8(d) shows the effect of correlation time $\tau_1$ on the PSD. It is clear that the interwell resonance becomes strong and the intrawell resonance becomes weak monotonously as $\tau_1$ increases. The reason is that the larger $\tau_1$ is, the narrower the bandwidth of colored noise is, then the more power spectrum passes through the system. These above results indicate that the additive colored noise plays a constructive role in occurrence of CR in the triple-well potential systems.

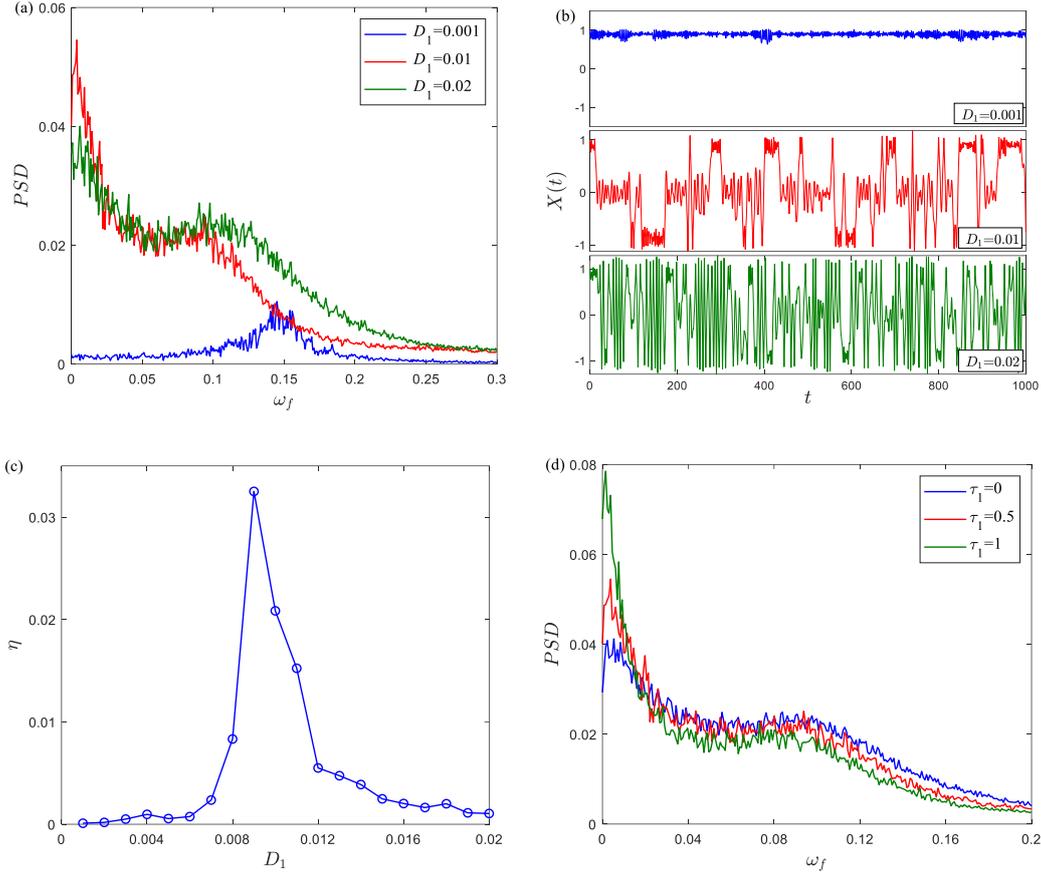

**FIG. 8.** (Color online) (a) The PSD as a function of frequency with different $D_1$; (b) the time series of system displacement with the different $D_1$; (c) the variation of the quality factor with $D_1$; (d) The PSD as a function of frequency with different $\tau_1$.

Then, the system (25) in case II is considered to explore the effects of the multiplicative colored noise and the cross-correlation between the additive and multiplicative colored noises on the CR. Figure 9 presents the variation of the PSD as a function of frequency $\omega_f$ with different noise intensity $D_2$ and colored noise cross-correlation $\lambda$. It is clear that the maximal value of PSD decreases monotonically with the increase of $D_2$ and increases monotonically as $|\lambda|$ increases. These indicate that the multiplicative colored noise and the colored noise cross-correlation cannot induce the CR. Thus, the additive and multiplicative colored noise present quite different effects for the features of system coherence.

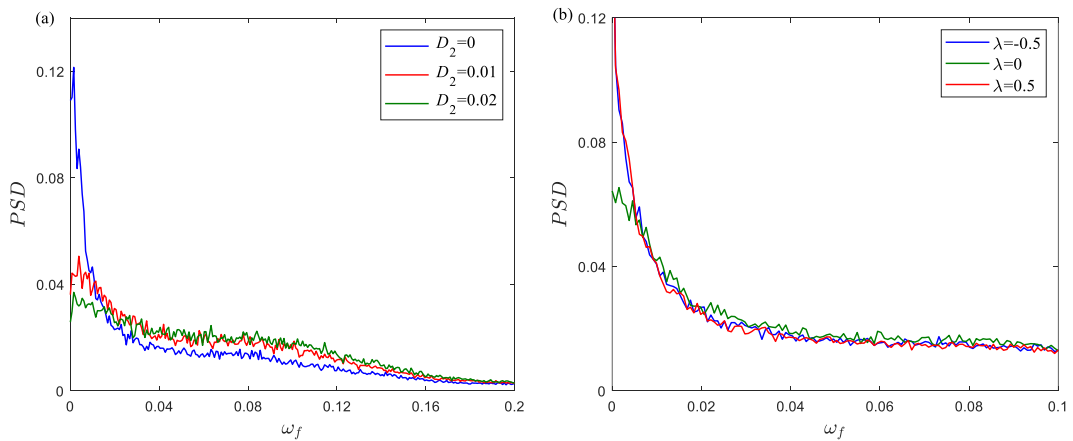

**FIG. 9.** (Color online) The PSD as a function of frequency with (a) different $D_2$ and (b) different $\lambda$.

## V. CONCLUSIONS

A stochastic averaging technique based on energy-dependent frequency is extended to obtain the analytical expression of SPD for the strongly nonlinear tri-stable stochastic systems. Then, two cases of strongly nonlinear tri-stable potential systems are presented to explore the effects of colored noise on stochastic dynamics and validate the effectiveness of the proposed method. Meanwhile, stochastic dynamical behaviors of the tri-stable systems induced by nonlinear stiffness coefficients are also discussed. The main conclusions can be drawn as follows. 1) The proposed method is well verified by numerical simulations, and has some significant advantages compared with the traditional stochastic averaging method, such as high accuracy, small limitations and easy application in tri-stable systems. 2) Noise intensity of colored noise can strengthen the transition among the three potential wells and enlarge the fluctuation range of stochastic response, which is opposite to the correlation time. 3) Nonlinear stiffness coefficients can change the topological structure of SPD and even induce stochastic P-bifurcation. 4) The additive and multiplicative colored noises display quite different effects on the features of CR. Choosing a moderate additive noise intensity can induce CR, but the multiplicative colored noise cannot. 5) Colored noise cross-correlation can break the symmetry of potential wells and induce phase transition. But it is not beneficial to the oscillatory regularity.

The extended stochastic averaging technique in this work plays a significant role in the theoretical study of tri-stable systems with random excitation which is not just limited to colored noise. What's more, based on energy-dependent frequency, the stochastic averaging technique can be also extended to multi-stable systems with potentials of more complicated shapes.

## ACKNOWLEDGEMENTS

This work is supported by the National Natural Science Foundation of China (Nos. 11772048, 11832005).

## REFERENCE

[1] M.S. Aswathy and S. Sarkar, Int. J. Mech. Sci. 153, 103 (2019).

[2] Q. F. He and M. F. Daqaq, J. Sound Vib. 333, 3479 (2014).

[3] M. Xu, Y. Wang, X. L. Jin, Z. L. Huang and T. X. Yu, Int. J. Non-Linear Mech. 52, 26 (2013).

[4] M. Xu and X. Y. Li, Int. J. Mech. Sci. 141, 206 (2018).

[5] S. M. Xiao and Y. F. Jin, Nonlinear Dyn. 90, 2069 (2017).

[6] W. Q. Zhu, Z. L. Huang and Y. Suzuki, Int. J. Nonlinear Mech. 36, 1235 (2001).

[7] H. Gang, T. Ditzinger, C.Z Ning, and H. Haken, Phys. Rev. Lett. 71, 807 (1993).

[8] Y. F. Jin, Z. M. Ma and S. M. Xiao, Chaos, Solitons and Fractals 103, 470 (2017).

[9] T. R. Bogatenko and V. V. Semenov, Phys. Lett. A 382, 2645 (2018).

[10] Y. F. Jin and P. F. Xu, IFAC. PapersOnLine 51,189 (2018).

[11] P. F. Xu, Y. F. Jin and S. M. Xiao, Chaos 27, 113109 (2017).

[12] R. L. Stratonovich, Discussion Papers 1, 30 (1967).

[13] R. Z. Khasminskii, Theory Probab. Appl. 11, 390 (1966).

[14] G. C. Papanicolaou and W. Kohler, Commun Math. Phys. 45, 217 (1975).

[15] G. Wainrib, Phys. Rev. E 84, 051113 (2011).

[16] P. D. Spanos, A. Sofi and M. D. Paola, J. Appl. Mech. 74, 315 (2007).

[17] J. R. Red-Horse and P. D. Spanos, Int. J. Non-Linear Mech. 27, 85 (1992).

[18] W. Q. Zhu and Y. K. Lin, J. Eng. Mech. 117, 1890 (1991).

[19] W. Q. Zhu, Z. L. Huang and Y. Q. Yang, J. Appl. Mech. 64, 975 (1997).

[20] Y. Xu, R. C. Gu, H. Q. Zhang, W. Xu and J. Q. Duan, Phys. Rev. E 83, 056215 (2011).

[21] T. Yang and Q. J. Cao, Mech. Syst. Signal Process 121, 745 (2019).

[22] I. S. M. Foukou, C. Nono, M. Siewe and C. Tchawoua, Commun. Nonlinear Sci. Numer. Simulat 56, 177 (2018).

[23] D. Liu, Y. Xu and J. L. Li, Chaos Solitons Fractals 104, 806 (2017).

[24] W. Q. Zhu, G. Q. Cai and R. C. Hu, Int. J. Dynamics and Control 1, 12 (2013).

[25] Y. F. Jin, S. M. Xiao and Y. X. Zhang, J. Stat. Mech. Theory E 2018, 123211 (2018).

[26] T. Yang and Q. J. Cao, Int. J. Mech. Sci. 156, 123 (2019).

[27] B. X. Fan and M.Xie, Phys. Rev. A 95, 023808 (2017).

[28] F. Bouthanoute, L. E. Arroum, Y. Boughaleb and M. Mazroui, M. J. Condensed Matter 9, 17 (2007).

[29] Y. Wang, Z. G. Ying and W. Q. Zhu, J. Sound Vib. 321, 976 (2009).

[30] W. Q. Zhu. Appl. Mech. Rev. 41, 189 (1988).

[31] Y. X. Zhang, Y. F. Jin, P. F. Xu and S. M. Xiao, Nonlinear Dyn. (2018).

[32] S. Arathi and S. Rajasekar, Phys. Scr. 84, 065011 (2011).

[33] L. Guan, Y. W. Fang, K. Z. Li, C. H. Zeng and F. Z. Yang, Phys. A 505, 716 (2018).